\documentclass[12pt,leqno]{article}
\usepackage[francais]{babel}
\usepackage[latin1]{inputenc}
\usepackage[T1]{fontenc}
\usepackage{amsmath,amsfonts,amssymb,amsthm,amscd}
\usepackage{bm}
\date{}

\newcommand\rg{\rightarrow}

\newcommand\N{\mathbb{N}}

\newcommand\p{\mathbb{P}}
\newcommand\R{\mathbb{R}}

\newcommand\CA{\mathcal{A}}
\newcommand\Cc{\mathcal{C}}

\newcommand\Ff{\mathcal{F}}
\newcommand\Gg{\mathcal{G}}

\newcommand\Bb{\mathcal{B}}

\newcommand\Oo{\mathcal{O}}

\newcommand\Rr{\mathcal{R}}

\newcommand\resp{\mathrm{resp}}

\def\adots{\mathinner{\mkern1mu\raise1pt\vbox{\kern7pt\hbox{.}}
\mkern2mu\raise4pt\hbox{.}
\mkern2mu\raise7pt\hbox{.}\mkern1mu}} 

\newtheorem{montheo}{\textsc{Th\'eor\`eme}}[section]
\newtheorem{monlem}{\textsc{Lemme}}[section]

\newtheorem{maprop}[monlem]{\textsc{Proposition}}

\begin{document}

\title{Etude du graphe divisoriel 6}

\author{Eric Saias}

%\date{juillet 2015}

\maketitle

\section{Introduction}

On appelle graphe divisoriel le graphe de la relation $\Rr_f$ définie sur l'ensemble des entiers positifs par

\begin{center}
$a$ $\Rr_f$ $ b$ si et seulement si $a$ divise $b$ ou $b$ divise $a$.
\end{center}

\noindent On s'intéresse plus spécifiquement ici aux propriétés asymptotiques quand $x$ tend vers $+\infty$ de la restriction $\Rr_f^x$ du graphe divisoriel aux entiers inférieurs ou égaux à $x$.

Sur l'ensemble des entiers inférieurs ou égaux à $x$, on définit également la relation $\Rr_g^x$ de la manière suivante où on note $[a,b]$ le plus petit commun multiple des entiers $a$ et $b$ :

\begin{center}
$a$ $\Rr_g^x$ $ b$ si et seulement si $[a,b]\le x$. 
\end{center}

\noindent On constate que si les entiers $a$ et $b$ sont en relation par $\Rr_f^x$, ils le sont également pour $\Rr_g^x$. Les graphes des relations $\Rr_f^x$ et $\Rr_g^x$ ont des propriétés communes. Cela justifie que l'on fasse quand c'est possible une étude simultanée de ces deux graphes.

On appelle chaîne d'entiers $\le x$ pour $\Rr_f$ (respectivement pour $\Rr_g^x$) de longueur $\ell$, tout $\ell$--uplet $(a_1,a_2,\ldots,a_\ell)$ d'entiers $\le x$ et deux à deux distincts tel que pour tout $i$ vérifiant $1\le i\le \ell$, on a $a_i$ $\Rr_f$ $a_{i+1}$ $(\resp.\ a_i$ $\Rr_g^x$ $a_{i+1})$. On notera $a_1-a_2-\cdots -a_\ell$ cette chaîne. On désigne par $f(x)$ $(\resp\ g(x))$ la longueur maximum d'une chaîne d'entiers $\le x$ pour $\Rr_f$ $(\resp\ \Rr_g^x)$. On note $\ell(\Cc)$ la longueur d'une chaîne $\Cc$.

Soient $\alpha$ et $\beta$ deux fonctions définies sur une partie $D$ de $\R^+$. On écrira $\alpha(x) \ll \beta(x)$ ou $\alpha(x)=O(\beta(x))$ quand il existe un réel positif $K$ tel que pour tout $x$ élément de $D$, on a $|\alpha(x)|\le K|\beta(x)|$. Quand on a simultanément $\alpha(x) \ll \beta(x)$ et $\beta(x)\ll \alpha(x)$, on écrira $\alpha(x)\asymp\beta(x)$.

A la suite des travaux de Pomerance \cite{Pom}, Pollington \cite{Pol} et Tenenbaum \cite{Ten2} sur l'asymptotique des fonctions $f(x)$ et $g(x)$, j'ai montré au théorème~2 de \cite{Sai2} que
\begin{equation}
f(x) \asymp g(x)\asymp x/\log x.
\end{equation}
Weingartner (Corollary 1 de \cite{Wei9}) a récemment donné une minoration explicite de $f(x)$ :
\begin{center}
\hspace{3cm}$f(x) > 0,76 \dfrac{x}{\log x}$, \hspace{2cm} pour $x$ suffisamment grand.
\end{center}

Au vu de (1.1), Erdös \cite{Erd} a posé la question de savoir comment ce résultat se généralise quand au lieu d'une seule chaîne, on en dispose de plusieurs. Cette question se subdivise en deux suivant que l'on demande aux chaînes d'être deux à deux disjointes ou non. De manière précise, désignons par $f(x,y)$ (respectivement $F(x,y)$) le cardinal maximum de l'ensemble des entiers représentés par l'union de $y$  chaînes quelconques (resp. deux à deux disjointes) d'entiers $\le x$ pour $\Rr_f$. Désignons également par $g(x,y)$ et $G(x,y)$ les analogues respectifs de $f(x,y)$ et $F(x,y)$ quand on remplace la relation $\Rr_f$ par la relation $\Rr_g^x$. Notons $\log^+(z) := \max (\log z,1)$. J'ai montré au théorème~1 de \cite{Sai4} que
\begin{equation}
\hspace{1cm}f(x,y) \asymp g(x,y) \asymp x \dfrac{\log^+y}{\log^+x},\hspace{3cm} \ \text{pour}\ x\ge y \ge1.
\end{equation}
Le principal objet de ce présent travail est de donner cette fois-ci l'ordre de grandeur des fonctions $F(x,y)$ et $G(x,y)$.

\vspace{2mm}

{\montheo {Il existe deux réels $K>c>0$ tels que l'on a pour tous $x\ge y \ge 1$}}
$$
c\,\frac{x}{\log^+ (x/y)}	 \le F(x,y) \le G(x,y) \le K \frac{x}{\log^+(x/y)}.
$$
On peut même demander que toutes les chaîne utilisées ont la même longueur, sans changer l'ordre de grandeur du nombre d'entiers représentés. Notons  $R(x,z)$ (respectivement $T(x,z)$) le nombre maximum d'entiers représentés par des chaînes deux à deux disjointes pour $\Rr_f$ $(\resp.\ \Rr_g^x)$ d'entiers $\le x$, et de longueur toutes égales à $z$.

\vspace{2mm}

{\montheo Il existe deux réels $K>c>0$ tels que pour tous entiers $x$ et $z$ vérifiant $1\le z\le f(x)$, on a}
$$
c\,\frac{x}{\log^+z}	 \le R(x,z) \le T(x,z) \le K \frac{x}{\log^+z}.
$$

Revenons à la fonction $F(x,y)$. La question sous--jacente à son introduction est de quantifier combien de nouveaux entiers on peut représenter quand on dispose de plusieurs chaînes disjointes et non d'une seule. Remarquons que l'on a trivialement $F(x,y) \le yf(x)$. Si cette majoration n'était pas très grossière, on pourrait  avoir  $\limsup\limits_{y\rg +\infty} \limsup\limits_{x\rg +\infty} (F(x,y)/f(x))=+\infty$ comme c'est le cas quand on remplace $F(x,y)$ par $f(x,y)$, d'après (1.2). Mais elle est de fait très grossière.  En effet le théorème~1 montre que l'on a en particulier $F(x,\sqrt{x})/f(x)\asymp 1$. Le nombre de nouveaux entiers représentés est donc très petit au regard du nombre de nouvelles chaînes dont on dispose.

\section{Les ensembles $\CA(x)$ et $\Bb(x)$ de Schinzel\break et Szekeres}

\setcounter{equation}{0}

Dans ce travail on utilisera la lettre $x$ pour désigner un élément générique de $\R^{+*}$, les lettres $m$ et $n$ pour désigner un élément générique de $\N^*$, et les lettres $p$ et $q$ pour désigner un nombre premier générique. Les lettres $i$, $j$, $k$ et $\ell$ serviront d'indice.

Pour tout $n\ge 2$, notons
$$
n=p_1p_2\cdots p_k\ \text{avec}\ p_1\ge p_2 \ge\cdots\ge p_k
$$
la décomposition de $n$ en facteurs premiers. On note $\Omega(n)=k$ et $P^-(n)=p_k$. On convient de plus que $\Omega(1)=0$ et $P^-(1)=+\infty$.

Notons à présent
\begin{equation}
S(n) := \left|
\begin{matrix}
&\max\limits_{1\le j\le k} \ p_1p_2\cdots p_{j-1}p_j^2 &\text{si} \ n\ge 2\\
&\text{et}\hfill\\
&1\hfill &\text{si}\ n=1
\end{matrix}
\right.
\end{equation}
la fonction introduite par Schinzel et Szekeres dans \cite{SchSze} en 1959, et
\begin{equation}
\CA(x):= \{n:S(n)\le x\}
\end{equation}
et
\begin{equation}
\Bb(x)\!:=\! \{2\le n\le x\! :\!\!\max\limits_{1\,\le\, j\,<\, k} p_1p_2\cdots p_{j-1}p_j^2\! \le\! x \!<\!p_1p_2\cdots p_{k-1}p_k^2 \}
\end{equation}
les ensembles qui lui sont directement associés. Le lien entre $\Bb(x)$ et $\CA(x)$ provient du fait que les entiers de $\Bb(x)$ sont les entiers $\le x$ n'appartenant pas  à $\CA(x)$, et minimaux pour la divisibilité (voir aussi les lemmes~2.2 et 2.4 ci--dessous). Remarquons dès à présent que l'on a
\begin{equation}
p\in \Bb(x) \Longleftrightarrow \sqrt{x} < p \le x. 
\end{equation}

Schinzel et Szekeres ont introduit ces objets pour répondre à une question d'Erdös \cite{SchSze}. Ruzsa \cite{Ruz} et Tenenbaum \cite{Ten1} et \cite{Ten2} ont montré qu'ils avaient plusieurs autres applications : aux entiers à diviseurs denses, au  petit crible d'Erdös et Ruzsa, aux nombres pratiques et au graphe divisoriel, (voir dans l'ordre chronologique \cite{Sai1}, \cite{Sai2}, \cite{Sai3},  \cite{Sai4}, \cite{Wei1}, \cite{Wei2}, \cite{Wei3}, \cite{Wei4}, \cite{Wei5}, \cite{Wei6}, \cite{Wei7}, \cite{TenWei}, \cite{Wei8}, \cite{Sai5} et \cite{Wei9}). On remarque que $S(n)$, $\CA(x)$ et $\Bb(x)$ sont souvent définis un peu différemment dans la littérature. Mais ce sont sous la forme (2.1), (2.2) et (2.3) que l'on travaillera avec ces objets dans  ce présent travail.

Précisons le lien avec les entiers à diviseurs denses en désignant par $1=d_1(n)<d_2(n) <\cdots <d_{\tau(n)}(n)=n$ la suite croissante de ses diviseurs. On dit que $n$ est à diviseurs $z$--denses quand
$$
\max_{1\le i <\tau(n)} \frac{d_{i+1}(n)}{d_i(n)} \le z.
$$
Le lien a été établi par Tenenbaum (lemme 2.2 de \cite{Ten1}) et est donné par la formule suivante :

{\monlem Soit $n\ge 2$. On a}
$$
S(n)=n\ \max_{1\le i<\tau(n)} \frac{d_{i+1}(n)}{d_i(n)}.
$$

L'étude des entiers à diviseurs denses est intéressante en soi (voir en particulier les travaux de Weingartner cités dans les références). Mais pour les questions relatives aux longues chaînes du graphe divisoriel (voir \cite{Ten2}, \cite{Sai2}, \cite{Sai4}, \cite{Sai5} et \cite{Wei9}) ce sont les propriétés et interactions des ensembles $\CA(x)$ et $\Bb(x)$ qui jouent le premier rôle. C'est donc encore le cas ici.

Les quatre résultats suivant donnent des propriétés fondamentales des ensembles $\CA(x)$ et $\Bb(x)$. Bien que les  trois premières soient déjà connues, nous pensons utile de regrouper ici leurs énoncés et de rappeler leurs preuves.

{\monlem Soit $1\le n\le x$ tel que $n\notin \CA(x)$. Il existe alors un unique élément de $\Bb(x)$ qui est aussi un diviseur de~$n$.}

\vspace{2mm}

\noindent \textbf{Démonstration} : Rappelons que l'on note $n=p_1p_2\cdots p_k$ avec $p_1\ge p_2 \ge \cdots \ge p_k$ la décomposition de l'entier $n\ge 2$ en facteurs premiers (on a bien ici $n\ge 2$ car $n\notin \CA(x)$ et $1\in \CA(x)$).

\vspace{2mm}

\noindent \textbf{Existence.}

L'entier $b:=p_1p_2\cdots p_j$ où $j$ désigne l'entier $i$ minimum tel que\break $p_1p_2\cdots p_{i-1}p_i^2>x$, est un élément de $\Bb(x)$ qui divise $n$.

\vspace{2mm}

\noindent \textbf{Unicité.}

Soit $b'=p_1^{\varepsilon_1}p_2^{\varepsilon_2}\cdots p_{\ell-1}^{\varepsilon_{\ell-1}}p_\ell$ avec $\ell\le k$ et $(\varepsilon_1,\varepsilon_2,\ldots ,\varepsilon_{\ell-1})\in \{0,1\}^{\ell-1}$, un élément de $\Bb(x)$ qui divise $n$. Supposons que $(\varepsilon_1,\varepsilon_2,\ldots,\varepsilon_{\ell-1}) \neq (1,1,\ldots,1)$. Comme $b'\in \Bb(x)$, on a alors
$$
x<p_1^{\varepsilon_1} p_2^{\varepsilon_2}\cdots p_{\ell-1}^{\varepsilon_{\ell-1}}p_\ell^2 \le p_1p_2\cdots p_\ell \le p_1p_2\cdots p_k \le x,
$$
ce qui est absurde. Donc
$$
b'=p_1p_2\cdots p_\ell.
$$
Notons $p_1p_2\cdots p_{j-2}p_{j-1}^2=1$ si $j=1$. En utilisant que $b$ et $b'$ sont éléments de $\Bb(x)$, on a $(p_1p_2\cdots p_{j-2})p_{j-1}^2 \le x < (p_1p_2\cdots p_{\ell-1})p_\ell^2$. Donc $\ell>j-1$ soit encore $\ell\ge j$. En inversant le rôle de $b$ et $b'$, on a aussi $j\ge \ell$. On a donc $\ell=j$ d'où $b'=b$.

{\monlem Soient $x>0$. On a pour tout couple $(b,b')$ d'entiers distincts de $\Bb (x)$,}
$$
[b,b']>x.
$$

\noindent \textbf{Démonstration} : Raisonnons par l'absurde en supposant que $[b,b']\le x$ pour deux entiers différents $b$ et $b'$ de $\Bb(x)$. Alors $b$ et $b'$ diviseraient le même entier $[b,b']\le x$. D'après le lemme 2.2, $[b,b']$ serait alors élément de $\CA(x)$, ce qui est absurde car tout diviseur d'un entier de $\CA(x)$ est aussi un élément de $\CA(x)$.

{\monlem On a pour tout $x>0$,}
$$
[x] = A(x)+\sum_{b\in \Bb(x)} \Big[\frac{x}{b}\Big].
$$

\noindent \textbf{Démonstration} : Cela découle immédiatement de la combinaison des\break lemmes~2.2 et 2.3.

{\monlem Notons
$$
e(t) := \left|
\begin{matrix}
&1 &\text{quand}\ [t]\ \text{est\ impair}\\
&0 &\text{quand}\ [t]\ \text{est\ pair}.\hfill
\end{matrix}
\right.
$$
On a alors pour tout $x>0$
$$
A(x)=\Big(\sum_{b\in \Bb(x)}e(x/b)\Big) - e(x).
$$ 
}

\noindent \textbf{Démonstration} : Rappelons que la formule
\begin{equation}
[x] = A(x) + \sum_{b\in \Bb(x)}\Big[\frac{x}{b}\Big]
\end{equation}
du lemme 2.4 résulte immédiatement de la combinaison des lemme 2.2 et 2.3. Comme les entiers de $\CA(x)$ sont inférieurs ou égaux à $x/2$, en criblant cette fois--ci uniquement les entiers $\le x/2$ par les éléments de $\Bb(x)$ et en utilisant de nouveau les lemmes 2.2 et 2.3, on a aussi
\begin{equation}
[x/2] = A(x) + \sum_{b\in \Bb(x)} \Big[\frac{x}{2b}\Big].
\end{equation}
Comme $e(t)=[t]-2[t/2]$, on obtient la formule annoncée en enlevant (2.5) au double de (2.6).

Rappelons ((2.4)) que les nombres premiers $p$ de $]\sqrt{x},x]$ sont éléments de $\Bb(x)$. En utilisant le théorème des nombres premiers et la formule $\log 2= \sum\limits_n
\frac{1}{2n-1}-\frac{1}{2n}$, on déduit du lemme 2.5 que

{\monlem On a pour $x \rg +\infty$}
$$
A(x) \ge (\log 2+o(1)) \frac{x}{\log x}.
$$

\noindent \textbf{Remarque 1.}  C'est nouveau : on sait maintenant que l'on peut minorer $A(x)$ par $(2/3)(x/\log x)$ qui est son ordre de grandeur, en utilisant d'une part le théorème des nombres premiers, et d'autre part de manière simple les propriétés fondamentales des ensembles $\CA(x)$ et $\Bb(x)$. Au contraire les preuves de  $A(x) \gg x /\log x$ qui étaient connues jusqu'à présent (voir d'une part la formule (1) de \cite{Sai2} qui résulte du travail \cite{Sai1}, et d'autre part la preuve de l'équivalent asymptotique du theorem 1 de \cite{Wei9}) sont longues et techniques.

\vspace{2mm}

\noindent \textbf{Remarque 2.} Revenons à présent au résultat
$$
f(x) \asymp x/\log x
$$ 
de la formule (1.1). Notons
$$
\CA(x,y) := \{n\in \CA(x) : P(n)\le y\}.
$$
J'ai montré au théorème 4.2 de \cite{Sai5}, que
\begin{equation}
f(x) \ge A(x/2)
\end{equation}
par la construction d'une chaîne qui utilise la propriété fondamentale des ensembles $\CA(x,y)$ donnée par la formule
$$
\CA(x,y) = \{1\} \bigsqcup\ \bigsqcup_{p\le \min (y,\sqrt{x})}p \CA(x/p,p),\qquad (x\ge 1,\ y\ge 1).
$$
En mettant bout à bout la minoration $f(x)\ge A(x/2)$ de la formule (2.7) et celle $A(x)\ge (0,69+o(1))x/\log x$ du lemme 2.6, on avance dans la compréhension des longues chaînes du graphe divisoriel. En effet on sait maintenant que les propriétés fondamentales de $\CA(x)$ (et de $\CA(x,y)$) mettent en lumière le fait que  l'ordre de grandeur $x/\log x$ de $f(x)$, provient de l'équivalence du théorème des nombres premiers : $\pi(x) \sim x/\log x$.

Généralisons $\CA(x)$ et $A(x)$ en posant $\CA(x,z,t)=\{n\le x : P^-(n)\ge z$ et $S(n)\le nt\}$ et $A(x,z,t)=|\CA(x,z,t)|$ (attention, on choisit ici pour ces notations l'inégalité large \og $P^-(n)\ge z$ \fg{} et non l'inégalité stricte \og $P(n)>z$ \fg{} comme dans \cite{Sai2} et \cite{Sai4}). On a  d'après les formules (1) et (3) et le lemme 5, tous de \cite{Sai2}, les résultats suivant

{\monlem On a
$$
\begin{array}{lll}
(i)\hspace{6cm} &\hspace{-4cm}A(x) \asymp \dfrac{x}{\log x}, &\hspace{3.5cm}(x\ge 2)\\
(ii) &\hspace{-4cm}A(x,z,t) \ll x \dfrac{\log^+(\min(x,t))}{\log^+ x \cdot \log^+ z}, &\hspace{1.7cm}(x\ge z\ge 1, t\ge 1)\\ \\
(iii) &\hspace{-4cm}B(x) \asymp x/\log x, &\hspace{3.5cm}(x\ge 2).
\end{array}
 $$}

\vspace{2mm}

\noindent \textbf{Remarque.} Il est à noter que Weingartner a donné une estimation précise de $A(x,z,t)$ (theorem 9 de \cite{Wei9}) qui fournit en particulier l'équivalent asymptotique de $A(x)$. Mais on n'a pas besoin ici de ses améliorations.

 {\monlem On a}
 $$
 \hspace{4cm}\sum_{b\in \Bb(x)}\frac{1}{b} = 1+O\Big(\frac{1}{\log x}\Big), \hspace{3.8cm} (x\ge 2).
 $$
 Cela résulte immédiatement de la combinaison du lemme 2.4 et des majorations de $A(x)$ et $B(x)$ fournies par les points \textit{(i)} et \textit{(iii)} du lemme 2.7. (voir encore dans \cite{Wei9}, le theorem~2 pour un résultat plus fin).
 
 On passe maintenant aux nouveaux résultats. Ce sont dans les lemmes 2.9 et 2.11 que réside la substantifique moëlle, respectivement de la minoration du théorème 2 et de la majoration du théorème 1. Le lemme 2.10 est quant à lui un résultat élémentaire d'analyse dont on se servira pour établir le lemme~2.11.
 
 {\monlem On a}
 $$
 \hspace{3cm}\sum_{\scriptstyle b\in \Bb(x)\atop \scriptstyle b>x/t} \frac{1}{b}\ll \frac{\log t}{\log x},\hspace{4cm}(x\ge 2,\ 2\le t\le x).
 $$
 
\vspace{2mm}

\noindent \textbf{Remarque.} Les calculs qui permettent d'établir cette majoration sont de même nature que ceux pour la majoration $B(x)\ll x/\log x$ (voir le lemme 6 de \cite{Sai2}) mais sensiblement plus longs. Voyons cela.

\vspace{2mm}

\noindent \textbf{Démonstration :} Cette majoration découle du lemme 2.8 quand\break $x^{1/4}<t\le x$. On supposera dorénavant
\begin{equation}
2\le t\le x^{1/4}.
\end{equation}

En utilisant (2.4) et (2.8), on a
$$
\sum_{\substack{ b\in \Bb(x)\\ \Omega(b)=1\\ b>x/t}} \frac{1}{b} = \sum_{x/t<p\le x}\frac{1}{p} \ll \frac{\log t}{\log x}.
$$
On a aussi
$$
\begin{matrix}
\displaystyle\sum_{\substack{b\in \Bb(x)\\ \Omega(b)=2\\ b>x/t}} \frac{1}{b} =\overset{}{\underset{\substack{qp\le x\\ q\ge p, q^2\le  x \\ qp^2>x, qpt>x}}{\displaystyle\sum\sum}} \frac{1}{qp}\\
\le \displaystyle\sum_{\frac{\sqrt{x}}{t}<p\le \sqrt{x}} \frac{1}{p} \sum_{x^{1/3}<q\le \sqrt{x}} \dfrac{1}{q}\asymp \dfrac{\log t}{\log x}.
\end{matrix}
$$

Supposons à présent $\Omega(b)\ge 3$ et écrivons  $b=aqp$ avec $P^-(a)\ge q\ge p$ et $a\not= 1$.

Les conditions $b\in \Bb(x)$ et $b>x/t$ se traduisent par $a\in \CA(\frac{x}{qp},q,qp)$, $aq^2\le x$, $aqp^2>x$ et $aqpt>x$. On a donc d'une part $p^3\le q^3 \le x$, d'où $qp \le x^{2/3}$, et d'autre part $q<p \, \min(p,t)$. Avec (2.8), l'inégalité $qp\le x^{2/3}$ entraîne également
$$
\frac{x}{pqt} \ge x^{1/12}.
$$
On a donc finalement
$$
\sum_{\scriptstyle  b\in \Bb(x), \Omega(b)\ge 3\atop\scriptstyle  b>x/t} \frac{1}{b} \le \sum_{p\le x^{1/3}} \frac{1}{p} \sum_{p\le q <p\min (p,t)} \frac{S(q,p)}{q}
$$
avec, en effectuant une sommation d'Abel et en utilisant le lemme 2.7\textit{(ii)},
$$
\begin{array}{ll}
&S(q,p) := \displaystyle\sum_{\scriptstyle a\in \CA(x/qp,\ q,\ qp)\atop\scriptstyle  a>\max(x/qp^2,\ x/qpt)} \dfrac{1}{a}\\\\
\le &\dfrac{A(x/qp,q,qp)}{x/qp} + \displaystyle\int_{\max(x/qp^2,\ x/qpt)} \dfrac{A(\theta,q,x/\theta)}{\theta^2}d\theta \\\\
\ll &\dfrac{1}{\log x} \Big(1+\displaystyle \int_{\max(x/qp^2,\ x/qpt)}^{x/qp} \dfrac{d\theta}{\theta}\Big)\\\\
 \asymp &\dfrac{\log(\min(p,t))}{\log x}.
\end{array}
$$
On a donc en utilisant (2.8) pour l'égalité
$$
\begin{array}{l}
\displaystyle\log x \sum_{\substack{ b\in \Bb(x)\\ \Omega(b) \ge 3\\ b>x/t}}\dfrac{1}{b} \ \ll \ \sum_{p\le x^{1/3}} \dfrac{1}{p} \sum_{p\le q< p\min(p,t)} \dfrac{\log(\min(p,t))}{q}\\\\
=\displaystyle \sum_{p\le t} \dfrac{\log p}{p} \sum_{p\le q<p^2}\dfrac{1}{q}
+\displaystyle \log t \sum_{t<p \le s^{1/3}} \dfrac{1}{p} \sum_{p\le q <pt}\dfrac{1}{q}\\
\asymp \log t.
\end{array}
$$
Cela conclut la preuve du lemme 2.9.

{\monlem Pour tous réels $s$ et $t$ vérifiant
\begin{equation}
t\ge 1\ \ et\ \ 0<s\le t
\end{equation}
on a
$$
\eta(s):= \eta(s,t) := s-1 -\frac{\log s}{1- \dfrac{\log s}{\log 9t}} \ge 0.
$$
}

\vspace{2mm}

\noindent\textbf{Démonstration :} On note $L=\log 9t$.

On a
$$
\eta'(s) = 1-\frac{1}{s(1-(\log s)/L)^2}
$$
et
$$
\eta''(s) = \dfrac{1-\dfrac{\log (e^2s)}{L}}{s^2(1-\dfrac{\log s}{L})^3}.
$$
Comme $9>e^2$, on a pour tous $s$ et $t$ vérifiant (2.9), $\eta''(s)>0$. On conclut en remarquant que $\eta(1)= \eta'(1)=0$.

{\monlem Il existe un réel $M\ge 2$ tel que pour tous entiers $x$ et $y$ vérifiant $x\ge M y \ge M$,  pour toute famille d'entiers $(m(b))_{b\in \Bb(x)}$ de $\N$ vérifiant
\begin{equation}
\forall b\in \Bb(x),\qquad 0\le m(b)\le x/b
\end{equation}
et
\begin{equation}
\sum_{b\in \Bb(x)} m(b) \le y,
\end{equation}
on a
$$
\sum_{b\in \Bb(x)} \frac{1}{b\, \log (9x/bm(b ))} \le \frac{1}{\log (9x/y)} - \frac{1}{5\log x}.
$$
}

\vspace{2mm}

\noindent\textbf{Démonstration :} Soit $(m(b))_{b\in \Bb(x)}$ une famille d'éléments de $\N$ vérifiant (2.10) et (2.11). On note
$$
\Bb^+(x) = \{b\in \Bb(x) : m(b) \ge 1\}
$$
et
$$
L= \log (9x/y).
$$
Pour tout $b\in\Bb^+(x)$, on a d'après (2.10) $\log(bm(b)/y)<L$ et
$$
\begin{array}{c}
\displaystyle \dfrac{1}{L}\Big(1+\dfrac{1}{L}\dfrac{\log(bm(b)/y)}{1-\dfrac{\log(bm(b)/y)}{L}}\Big)\\
= \dfrac{1}{L}\Big(1+ \dfrac{\log(bm(b)/y)}{\log(9x/bm(b))}\Big) = \dfrac{1}{\log(9x/bm(b))}
\end{array}
$$
D'où, en utilisant (2.11) pour l'inégalité,
$$
\begin{array}{c}
\displaystyle \dfrac{1}{L} - \sum_{b\in \Bb^+(x)} \dfrac{1}{b\, \log(9x/bm(b))}\\
= \displaystyle \dfrac{1}{L}\Big[ 1- \sum_{b\in \Bb^+(x)} \dfrac{1}{b} \Big( 1+ \dfrac{1}{L} \dfrac{\log(bm(b)/y)}{1-\dfrac{\log(bm(b)/y)}{L}} \Big)\Big]\\
\ge \displaystyle \dfrac{1}{L} \Big(1-\sum_{b\in \Bb^+(x)} \dfrac{1}{b}\Big) + \dfrac{1}{L^2} \Big[ -1+ \sum_{b\in \Bb^+(x)}\dfrac{1}{b} \Big(\dfrac{bm(b)}{y} - \dfrac{\log(bm(b)/y)}{1- \dfrac{\log(bm(b)(y)}{L}} \Big)\Big]
\end{array}
$$
soit finalement
\begin{equation}
\begin{array}{l}
\displaystyle \dfrac{1}{\log(9x/y)} - \sum_{b\in\Bb^+(x)} \dfrac{1}{b\, \log (9x/bm(b))}\\
\ge \displaystyle \dfrac{1}{L} \Big(1- \dfrac{1}{L}\Big) \Big(1- \sum_{b\in \Bb^+(x)} \dfrac{1}{b}\Big) + \dfrac{1}{L^2} \sum_{b\in \Bb^+(x)} \dfrac{\eta(bm(b)/y, x/y)}{b}.
\end{array}
\end{equation}

\vspace{2mm}

\noindent\textbf{Remarque 1.} 

Cette dernière inégalité constitue le point central de la preuve.

\vspace{2mm}

\noindent\textbf{Remarque 2.} 

Quand $s=1$, $\eta(s)$ prend sa valeur minimum qui est nulle. En supposant que la fonction $m$ soit plus généralement à valeurs dans $\R^+$, le second terme à droite du signe $\ge$ de la dernière inégalité est nul avec le choix $m(b)=y/b$. Alors, même dans le cas le plus favorable qui est $1-\sum\limits_{b\in \Bb^+(x)} \dfrac{1}{b} \asymp \dfrac{1}{\log x}$, l'expression à droite du signe $\ge$ serait $\asymp \dfrac{1}{L\log x}$, ce qui est insuffisant pour conclure. On verra donc que l'on fait  dans la preuve un usage crucial de l'hypothèse que $m$ est à valeurs entières.

Reprenons le fil de la preuve du lemme 2.11. La combinaison de l'inégalité (2.12) et des lemmes 2.8 et 2.10 permet de conclure dans le cas où $\sum\limits_{b\in \Bb(x)\atop m(b)=0} \dfrac{1}{b}>\dfrac{L}{4\, \log x}$ et $x/y$ est suffisamment grand. Si au contraire $\sum\limits_{b\in\Bb(x)\atop m(b)=0}\dfrac{1}{b}\le \dfrac{L}{4\, \log x}$, on~a alors en supposant toujours que $x/y$ est suffisamment grand,
$$
\begin{array}{ll}
 &\displaystyle\sum_{b\in \Bb^+(x)} \dfrac{\eta(bm(b)/y,x/y)}{b} \ge \sum_{b\in \Bb^+(x)\atop b>\sqrt{xy}} \dfrac{\eta(bm(b)/y)/y,x/y)}{b}\\
&\displaystyle\ge \dfrac{12}{13}\sqrt{\dfrac{x}{y}} \sum_{b\in \Bb^+(x)\atop b>\sqrt{xy}} \dfrac{1}{b}\\
&\displaystyle\ge \dfrac{12}{13}\sqrt{\dfrac{x}{y}}\Big( \sum_{b\in \Bb(x)\atop b>\sqrt{xy}} \dfrac{1}{b} - \dfrac{L}{4\,\log x}\Big)\\
&\displaystyle\ge \dfrac{12}{13}\sqrt{\dfrac{x}{y}}\Big( \sum_{\sqrt{xy}<p\le x} \dfrac{1}{p} - \dfrac{L}{4\,\log x}\Big)\\
&\ge \dfrac{1}{13} \sqrt{\dfrac{x}{y} }\dfrac{L}{\log x}.
\end{array}
$$
On en déduit en utilisant (2.12) et le lemme 2.8 que, toujours pour $x/y$ suffisamment grand, on a
$$
\begin{array}{l}
\displaystyle \dfrac{1}{\log(9x/y)} - \sum_{b\in \Bb(x)} \dfrac{1}{b \log(9x/bm(b))}\\
\displaystyle \ge \dfrac{1}{13 \log x} \dfrac{\sqrt{x/y}}{\log (9x/y)} - \dfrac{1}{\log x} \ge \dfrac{1}{5 \log x}.
\end{array}
$$

\section{Composantes}

Soit $\Cc= n_1-n_2\cdots - n_\ell$ une chaîne d'entiers $\le x$ pour $\Rr_g^x$. Rappelons (voir le cas particulier $t=1$ du chapitre 6 de \cite{Sai2}) que l'on appelle composante de $\Cc$ toute chaîne extraite de $\Cc$ de la forme $\Oo:= n_i-n_{i+1}-\cdots -n_j$ avec $1\le i\le j \le \ell$ et vérifiant
$$
\begin{array}{l}
n_\alpha\notin \CA(x) \ \ \text{pour \ }i\le \alpha\le j\\
i=1\ \text{ou}\ n_{i-1}\in \CA(x)\\
j=\ell \ \text{ou}\ n_{j+1}\in \CA(x).
\end{array}
$$
Le cas particulier $t=1$ du lemme 11 de \cite{Sai2} correspond au résultat suivant.

{\monlem Soient $x\ge 2$ et $\Cc$ une chaîne d'entiers $\le x$ pour $\Rr_g^x$. Soit $\Oo$ une composante de $\Cc$. Il existe  alors un unique entier $b$ de $\Bb(x)$ qui divise tous les entiers de $\Oo$. On le note $b(\Oo)$.}

\vspace{2mm}

\noindent\textbf{Démonstration :} Comme il est indiqué dans \cite{Sai2}, ce résultat découle facilement des lemmes 2.2 et 2.3.

\vspace{2mm}

\noindent\textbf{Remarque.} Une chaîne d'entiers $\le x$ au sens de $\Rr_f$ est aussi une chaîne  d'entiers $\le x$ au sens de $\Rr_g^x$. Donc la notion de composante et le lemme 3.1 sont valables aussi pour des chaînes au sens de~$\Rr_f^x$.

\section{Minoration de $R(x,z)$}

\noindent\textbf{a) Inéquation fonctionnelle pour $R(x,z)$}

Notons $f_a(x)$ la longueur maximum d'une chaîne d'entiers de $\CA(x)$.

{\maprop  Pour tout couple d'entiers positifs $x$ et $z$ tels que $1\le z\le f_a(x)$, on a
$$
R(x,z) \ge \sum_{b\in \Bb(x)} R\Big(\frac{x}{b},z\Big)\ \ \  + \frac{f_a(x)}{2}.
$$}

\vspace{2mm}

\noindent\textbf{Démonstration :} Pour tout entier $b$ fixé de $\Bb(x)$, il existe $R(x/b,z)$ (éventuellement nul) chaînes d'entiers $\le x$ deux à deux disjointes, toutes de longueur $z$, et de la forme $b\,\Cc(x/b)$ où $\Cc(x/b)$ désigne une chaîne d'entiers $\le x/b$. De plus quand on fait maintenant varier $b$ dans $\Bb(x)$, on déduit du lemme 2.3 que toutes ces $\sum\limits_{b\in \Bb(x)}R(\dfrac{x}{b},z)$ chaînes sont deux à deux disjointes.

Par ailleurs notons $\Cc_a(x)$ une chaîne d'entiers de $\CA(x)$ et de longueur $f_a(x)$. On effectue la division euclidienne de $f_a(x)$ par $z$
$$
f_a(x) = dz+r \quad \text{avec} \quad 0\le r<z
$$
et on extrait de $\Cc_a(x)$, $d$ sous--chaînes deux à deux disjointes et toutes de longueur $z$. D'après l'hypothèse $f_a(x)\ge z$, le nombre d'entiers représentés par ces $d$ chaînes est égal à $dz>f_a(x)/2$. De plus  un multiple d'un élément de $\Bb(x)$ ne peut pas appartenir à $\CA(x)$. On en déduit que ces $d$ chaînes sont disjointes des $\sum\limits_{b\in \Bb(x)}R(x/b,z)$ chaînes construites précédemment. Cela conclut la preuve de la proposition 4.1.

\vspace{2mm}

\noindent\textbf{b) Minoration du second membre}

{\monlem On a 
$$
\hspace{3cm}f_a(x)\gg x/\log x, \hspace{6cm}(x\ge 2).
$$}

\vspace{2mm}

\noindent\textbf{Démonstration :} Cela résulte du théorème 5.1 de \cite{Sai5}.

\vspace{2mm}

\noindent\textbf{c) Minoration de $R(x,z)$}

L'objet de ce chapitre consiste à établir la minoration
\setcounter{equation}{0}
\begin{equation}
 \hspace{2.5cm}R(x,z) \gg \frac{x}{\log ^+ z},  \hspace{3.5cm}(1\le z\le f(x))
\end{equation}
du théorème 2.

Il découle des lemmes 2.8 et 2.9, et (1.1), qu'il existe trois gros réels $K_1$, $K_2$ et $z_0$ tels que
\begin{eqnarray}
\displaystyle \sum_{b\in \Bb(x)} \frac{1}{b} \ge 1-\frac{K_1}{\log x}, \hspace{4cm}(x\ge 2)\\
\displaystyle \sum_{b\in \Bb(x)\atop b> x/t} \frac{1}{b} \le K_2 \frac{\log t}{\log x},  \hspace{3cm}(x\ge 2,\ t\ge 2)\\
\hspace{2cm}
K_1+K_2\log(z\,\log^2 z) \le 2K_2\, \log z,  \hspace{3.5cm}(z\ge z_0)
\end{eqnarray}
et
\begin{equation}
[\sqrt{x} \ge z \, \log^2 z\ \text{et}\ z\ge z_0] \ \text{entraînent}\ f(\sqrt{x})\ge z.
\end{equation}
On distingue maintenant quatre cas pour établir (4.1).

\hspace{2mm}

\noindent\textbf{1\up{er} cas :} $z<z_0$ et $x<2^z$.

Comme $f(x)\ge z$, on peut extraire d'une chaîne d'entiers $\le x$ de longueur $f(x)$ une chaîne de longueur $z$. Cela permet de conclure dans ce cas.

\hspace{2mm}

\noindent\textbf{2\up{ième} cas :} $z<z_0$ et $x\ge 2^z$.

A tout nombre impair $m$ vérifiant $1\le m\le x/2^z$, on associe la chaîne de longueur $z$ :
$$
2m - 4m -8 m -\cdots - 2^z m.
$$
Ces chaînes étant deux à deux disjointes, on en déduit que
$$
R(x,z) \ge \dfrac{z}{2^{z+1}} x \asymp \dfrac{x}{\log z}.
$$

\hspace{2mm}

\noindent\textbf{3\up{ième} cas :} $f_a(x) <z\le f(x)$ et $z\ge z_o$.

On extrait d'une chaîne d'entiers $\le x$ de longueur $f(x)$ une sous chaîne de longueur $z$. Avec le lemme 4.2, on a
$$
R(x,z) \ge z >f_a(x) \gg x/\log x \asymp x/\log z.
$$

\hspace{2mm}

\noindent\textbf{4\up{ième} cas :} $z_0 \le z \le f_a(x)$.

En combinant la proposition 4.1 et le lemme 4.2, on peut choisir un petit réel $c_1>0$ tel que
\begin{equation}
R(x,z) \ge \sum_{b\in \Bb(x)} R\Big(\dfrac{x}{b},z\Big)\ \ \  + c_1 \frac{x}{\log x}.
\end{equation}

D'après les trois premiers cas, il existe un réel $c>0$ tel que
\begin{equation}
c\le c_1/2K_2
\end{equation}
et
\begin{equation}
\begin{matrix}
\hspace{1cm}R(x,z) \ge cx/\log^+ z \hspace{3cm} (1\le z<z_0\ \text{ou}\ f_a(x)<z).
\end{matrix}
\end{equation}

Montrons alors par récurrence sur $k\ge 0$ que l'on a
\begin{equation}
\left| \begin{matrix}
R(x,z) \ge c \dfrac{x}{\log^+z}, \ \text{dans\ la\ région}\\ \\
(H_k)\ 2\le x <2^{2^k} \ \text{et}\ 1\le z\le f(x).
\end{matrix}
\right.
\end{equation}

Le cas où $k=0$ est vide. Supposons à présent (4.9) vérifiée, ainsi que la condition $(H_{k+1})$. D'après (4.8), on peut supposer que
$$
z_0 \le z \le f_a(x).
$$
En utilisant successivement l'hypothèse de récurrence (4.9) avec (4.5), (4.3), (4.2) et (4.4), on a
$$
\begin{matrix}
\displaystyle \frac{1}{c} \sum_{\scriptstyle b\in \Bb(x)\atop \scriptstyle b\le \frac{\scriptstyle x}{\scriptstyle z\log ^2z}} R(\frac{x}{b},z) &\ge\displaystyle \frac{x}{\log^+z} \sum_{\scriptstyle b\in \Bb(x) \atop \scriptstyle b\le \frac{\scriptstyle x}{\scriptstyle z\,\log^2z}} \frac{1}{b}\hfill\\
&\ge \displaystyle \frac{x}{\log^+z}\Big(\sum_{b\in \Bb(x)}\frac{1}{b}-\ \ \  \frac{K_2 \log (z\log^2 z)}{\log x}\Big)\hfill\\
&\ge \displaystyle \frac{x}{\log^+z}\Big(1- \frac{1}{\log x} (K_1+K_2 \log(z\log^2 z))\Big)\\ \\
& \ge \dfrac{x}{\log^+z} - 2K_2 \dfrac{x}{\log x}.\hfill
\end{matrix}
$$
D'où avec (4.6)
$$
R(x,z) \ge c\,\frac{x}{\log^+z} +(c_1-2K_2 c) \frac{x}{\log x}\ge c\ \dfrac{x}{\log^+z}
$$
d'après (4.7). Cela achève la partie itérative de la récurrence, et partant la preuve de (4.1).

\section{Majoration de $G(x,y)$}

Dans tout ce chapitre, les chaînes d'entiers $\ge x$ seront implicitement des chaînes d'entiers au sens de $\Rr_g^x$.

La preuve de la majoration de $G(x,y)$ suit les mêmes étapes que celle de la minoration de $R(x)$ : inéquation fonctionnelle, majoration du second membre qui permet d'initialiser la récurrence, et enfin la partie itérative de la récurrence.

\vspace{2mm}

\noindent\textbf{a) Inéquation fonctionnelle pour $G(x,y)$}

{\maprop Pour tous entiers $x$ et $y$ tels que $x\ge 2y\ge 2$, on a
$$
G(x,y) \le \max_{m} \sum_{b\in \Bb(x)} G\Big(\frac{x}{b},m(b)\Big) + \max_{\Gg(x)} \sum_{\scriptstyle \Cc\in \Gg(x)\atop \scriptstyle \Cc\cap \CA(x)\not= \emptyset}\ell(\Cc)
$$
où le premier maximum porte sur les familles d'entiers $(m(b))_{b\in \Bb(x)}$ de $\N$ qui vérifient
\setcounter{equation}{0}
\begin{equation}
\forall\ b\in \Bb(x),\quad 0\le m(b)\le x/b
\end{equation}
et
\begin{equation}
\sum_{b\in\Bb(x)} m(b) \le y,
\end{equation}
et le second maximum porte sur toutes les partitions $\Gg(x)$ de l'ensemble des entiers $\le x$ en chaînes d'entiers $\le x$.} 

\vspace{2mm}

\noindent\textbf{Démonstration :} Soit $\Ff$ une famille de $y$ chaînes d'entiers $\le x$ deux à deux disjointes, et dont l'union est de cardinal maximum $G(x,y)$. Toutes celles qui ne contiennent aucun entier de $\CA(x)$ sont des composantes d'elles--mêmes. Par  le lemme 3.1, elles sont donc de la forme $b\,\Cc(x/b)$ où $b\in \Bb(x)$ et $\Cc(x/b)$ est une chaîne d'entiers $\le x/b$. Le nombre d'entiers représentés par toutes ces chaînes est donc inférieur ou égal à
$$
\max_m \sum_{b\in \Bb(x)} G(x/b,m(b))
$$
où le maximum est défini dans la proposition.

Par ailleurs, on complète l'union disjointe des chaînes de la famille $\Ff$ qui rencontrent $\CA(x)$ en une partition $\Gg(x)$ des entiers $\le x$, en rajoutant des chaînes--singletons. Cela achève de prouver la proposition 5.1.
 
\vspace{2mm}

\noindent\textbf{b) Majoration du second membre}

{\monlem On a pour $x\ge 2$,
$$
\max_{\Gg(x)} \sum_{\scriptstyle \Cc\in \Gg(x)\atop \scriptstyle \Cc\cap \CA(x)\not=\emptyset} \ell(\Cc) \ll x/\log x
$$
où le maximum porte sur toutes les partitions $\Gg(x)$ de l'ensemble des entiers $\le x$ en chaînes d'entiers $\le x$.}

\vspace{2mm}

\noindent\textbf{Remarque.} A quelques détails près, c'est en réalité ce lemme qui est démontré quand on établit la majoration $g(x)\ll x/\log x$ au chapitre 7 de \cite{Sai2}. Voyons cela.

\vspace{2mm}

\noindent\textbf{Démonstration :} Soit $\Cc=n_1-n_2-\cdots -n_\ell$ une chaîne d'entiers.  On note $\chi(\Cc)$ l'ensemble des composantes de $\Cc$ qui ne commencent pas en $n_1$. Soit $\Gg(x)$ une partition des entiers $\le x$ en chaînes d'entiers $\le x$. 

Pour majorer $g(x)$ dans \cite{Sai2}, on a utilisé que
\begin{equation}
card \bigcup_{\Oo\in \chi(\Cc)}
\Oo \le \sum_{a\in \CA(x)}(S(x,a)+T(x,a))
\end{equation}
(voir le milieu de la page 239 de \cite{Sai2} pour cette inégalité et la définition de $S(x,a)$ et $T(x,a)$).
Pour énoncer l'analogue de la majoration (5.3) dont on a besoin ici, on est amené à changer certaines des notations de \cite{Sai2}.

Soit $\Gg(x)$ une partition des entiers $\le x$ en chaînes d'entiers inférieurs ou égaux à $x$. Pour $\Oo\in \bigsqcup\limits_{\Cc\in \Gg(x)} \chi(\Cc)$, on garde les mêmes notations que dans \cite{Sai2} pour $e(\Oo)$, $b(\Oo)$, $a(\Oo)$ et $r(\Oo)$. En revanche pour tout $a\in \CA(x)$, on note
$$
\begin{array}{ll}
\chi(\Gg(x),a) &= \Big\{ \Oo \in \bigcup\limits_{\Cc\in \Gg(x)} \chi(\Cc) : a(\Oo)=a\Big\}\\
\Rr(\Gg(x),a) &=\Big\{ r(\Oo) : \Oo \in \chi (\Gg(x),a)\Big\}
\end{array}
$$
et
$$
k(a,r) = card\Big\{ \Oo \in \bigsqcup_{\Cc \in \Gg(x)} \chi(\Cc) : a=a(\Oo)\ \text{et}\ r=r(\Oo)\Big\}.
$$
On note enfin
$$
S'(x,a) = \sum_{\scriptstyle r\in \Rr(\Gg(x),a)\atop\scriptstyle  k(a,r)>\log^3(x/a)\ \text{ou}\ r>(x/a)^{2/3}} \frac{x}{ar}
$$
et
$$
T'(x,a) = \sum_{\scriptstyle r\in \Rr(\Gg(x),a)\atop\scriptstyle  r\le (x/a)^{2/3}} g(x/ar,\log^3(x/a)).
$$

De manière analogue à (5.3), on a
\begin{equation}
card \bigsqcup_{\Oo\in \bigsqcup\limits_{\Cc\in \Gg(x)}\chi(\Cc)} \Oo \le \sum_{a\in \CA(x)} (S'(x,a)+T'(x,a)).
\end{equation}
Les inégalités suivantes correspondent à celles des lemmes 14 et 15 de \cite{Sai2}.

\vspace{2mm}

\noindent\textbf{Lemme 14'.} \textit{Soient $\Gg(x)$ une partition des entiers $\le x$ en chaînes d'entiers inférieurs ou égaux à $x$, et $a\in \CA(x)$. Pour tout élément $r$ de de $\Rr(\Gg(x),a)$, on a}
$$
\max \Big(\frac{x}{aP^-(a)},\Big(\frac{x}{a}\Big)^{\Omega(r)/(\Omega(r)+1)}\Big) < r \le \frac{x}{a}.
$$
et
\vspace{2mm}

\noindent\textbf{Lemme 15'.} \textit{Soient $\Gg(x)$ une partition des entiers $\le x$ en chaînes d'entiers $\le x$ et $a\in \CA(x)$. On a}
$$
card\ \Rr(\Gg(x),a) \le card\ \chi(\Gg(x),a) <\min (P^-(a),\sqrt{x/a}).
$$

En suivant la démarche de \cite{Sai2} pour établir la majoration
$$
\sum_{a\in \CA(x)} S(x,a) + T(x,a) \ll x/\log x,
$$
en gardant les lemmes 7, 13 et 16 de \cite{Sai2} inchangés, et en remplaçant les lemmes 14 et 15, $\chi(\Cc,a)$ et $\Rr(\Cc,a)$ par respectivement les lemmes 14' et 15', $\chi(\Gg(x),a)$ et $\Rr(\Gg(x),a)$, on obtient  de même
$$
\sum_{a\in \CA(x)} S'(x,a) + T'(x,a) \ll x/\log x.
$$
En combinant avec (5.4), on a donc montré que
\begin{equation}
card \bigsqcup_{\scriptstyle \Oo\in \bigsqcup \chi(\Cc)\atop \scriptstyle \Cc\in \Gg(x)}\Oo \ll x\log x.
\end{equation}
On a donc aussi
\begin{equation}
card \bigsqcup_{\scriptstyle \Oo\in \bigsqcup \chi(\Cc)\atop\scriptstyle  \Cc\in \Gg(x)}\Oo = card \bigsqcup_{\scriptstyle \Oo\in \bigsqcup\chi(\Cc)\atop \scriptstyle \Cc \in inv \Gg(x)} \Oo \ll x/\log x.
\end{equation}

Soit maintenant $\Gg(x)$ une partition des entiers $\le x$ en chaînes d'entiers $\le x$, $\Cc=n_1-n_2-\cdots - n_\ell$ une chaîne de $\Gg(x)$ et $\Oo$ une composante de $\Cc$. On a alors une trichotomie. Soit $\Oo$ ne commence pas en $n_1$ ; soit $\Oo$ ne finit pas en $n_1$ ; soit $\Oo$ n'est formée que d'un seul entier qui appartient à $\Bb(x)$. On a donc
$$
\begin{matrix}
\displaystyle\sum_{\scriptstyle \Cc\in \Gg(x)\atop\scriptstyle  \Cc\cap\CA(x)\not=\emptyset} \ell(\Cc) &\displaystyle\le A(x) + B(x) + 2 \displaystyle\sum_{\scriptstyle \Oo\in \bigsqcup \chi(\Cc)\atop \scriptstyle\Cc\in \Gg(x)}\ell(\Oo) \\
&\asymp \dfrac{x}{\log x}\hfill
\end{matrix}
$$
d'après (5.5) et les points (i) et (iii) du lemme 2.7. Cela conclut.

\vspace{2mm}

\noindent\textbf{Remarque.} On a donc ici été amené à réutiliser le lemme 16 de \cite{Sai2}. Même si cela n'est pas utile ici, signalons que l'on peut légèrement améliorer ce lemme par le résultat suivant

\vspace{2mm}

\noindent\textbf{Lemme 16*.} \textit{Il existe un réel $C$ tel que pour toute partie finie non vide $\p$ de l'ensemble des nombres premiers, on a}
$$
\sum_{p\in \p} \frac{1}{p} \le C \frac{\log(1+\#\p\cdot \frac{\log(\min \p)}{\min \p})}{\log(\min \p)}.
$$

\vspace{2mm}

\noindent\textbf{c) Majoration de $G(x,y)$}

L'objet de ce chapitre est d'établir la majoration
\begin{equation}
\hspace{1.4cm}G(x,y) \ll \frac{x}{\log(x/y)}, \hspace{5cm} (x\ge 2y\ge 2)
\end{equation}
du théorème 1.

Soit $k_o$ un entier tel que
\begin{equation}
2^{2^{k_o}}\ge M
\end{equation}
où $M$ est le réel du Lemme 2.11. Soit $K_3$ le réel implicite correspondant au symbole $\ll$ dans le lemme 5.2. On choisit alors un réel $K$ vérifiant
\begin{equation}
K\ge 5K_3 
\end{equation}
et
\begin{equation}
\hspace{1cm}G(x,y) \le K\frac{x}{\log(9x/y)},\hspace{3.3cm} (1\le y \le x \le 2^{2^{k_o}}y).
\end{equation}
Montrons à présent par récurrence sur $k\ge k_o$ que
$$
\hspace{1cm} G(x,y) \le K \frac{x}{\log(9x/y)}, \hspace{3cm} (1\le y \le x \le 2^{2^k}y).\leqno(H_k)
$$

L'initialisation provient de (5.10). De plus le cas où $x/y< M$ provient de la combinaison des inégalités (5.10) et (5.8). On peut donc supposer dorénavant
\begin{equation}
x\ge My.
\end{equation}

Supposons maintenant $(H_k)$ vérifiée pour un certain $k\ge k_o$. Dans les inégalités ci--dessous, le maximum porte sur les familles d'entiers $(m(b))_{b\in \Bb(x)}$ de $\N$ vérifiant (5.1) et (5.2). En utilisant successivement la proposition 5.1 et le lemme 5.2, l'hypothèse de récurrence, le lemme 2.11 avec (5.11) et enfin (5.9), on obtient pour $1\le y \le x\le 2^{2^{k+1}}y$,
$$
\begin{array}{rl}
G(x,y) &\le K_3 \displaystyle \frac{x}{\log x} + \max\limits_m \sum_{b\in \Bb(x)} G(\frac{x}{b},m(b))\\
&\displaystyle\le x\Big[\frac{K_3}{\log x} + K \max\limits_m \sum_{b\in \Bb(x)} \frac{1}{b\,\log(9x/bm(b))}\Big]\\
&\displaystyle\le x \Big[\frac{K}{\log (9x/y)}	+ (K_3 -\frac{K}{5})\frac{1}{\log x}\Big]\\ \\
&\le \dfrac{Kx}{\log (9x/y)}.
\end{array}
$$
Cela achève l'étape itérative de la récurrence, et partant la preuve de (5.7).

\section{Fin de la preuve des théorèmes 1 et 2}

\setcounter{equation}{0}

On a prouvé aux chapitres 4 et 5 respectivement la minoration du théorème 2 et la majoration du théorème 1.

Par ailleurs, on a pour $x\ge y\ge 1$
\begin{eqnarray}
&&F(x,y) \ge f(x)\\
&&F(x,y) \ge R(x,\lceil x/y\rceil)
\end{eqnarray}
et
\begin{equation}
T(x,z) \le G(x,x/z).
\end{equation}
On observe alors que la minoration du théorème 1 découle de la combinaison de (6.1) et (1.1) quand $y\le \sqrt{x}$, de celle de (6.2) et de la minoration du théorème 2 quand $y>\sqrt{x}$ et $x$ est suffisamment grand, et est finalement banale quand $x$ est borné. De même la majoration du théorème 2 résulte de la combinaison de (6.3) avec la majoration du théorème 1. Cela achève donc la preuve des théorèmes 1 et 2.

\vspace{1cm}

\begin{center}
Remerciement
\end{center}

\vspace{2mm}

Je remercie Pierre Mazet de m'avoir donné l'idée d'étudier la fonction $R(x,z)$.

\vskip4mm

\begin{tabular}{ll}

 &\hspace{6.7cm}Eric Saias \\

 &\hspace{6.7cm}Sorbonne Université\\

&\hspace{6.7cm}LPSM\\

&\hspace{6.7cm}4, place Jussieu\\

&\hspace{6.7cm}75252 Paris Cedex 05 (France)\\

\vspace{2mm}

 &\hspace{6.7cm}\textsf{eric.saias@upmc.fr}
\end{tabular}

\end{document}